\documentclass[12pt]{article}
\usepackage[cp850]{inputenc}
\usepackage[T1]{fontenc}\usepackage{amsfonts}
\usepackage{graphicx}
\usepackage{amsmath}

\title{The mean perimeter  of some random plane convex sets generated by  a Brownian motion}
\author{Philippe Biane\thanks{CNRS, Laboratoire d'Informatique, Institut Gaspard Monge, Universit\'e Paris-Est, 5 Boulevard Descartes, Champs-sur-Marne F-77454, Marne-la-Vallee Cedex 2 \texttt{Philippe.Biane@univ-mlv.fr}}\\ G\'erard Letac\thanks{Laboratoire de Statistique et Probabilit\'es,
Universit\'e Paul Sabatier, 31062, Toulouse, France, \texttt{letac@cict.fr}}}
\date{}

\def\<{\langle}
\def\>{\rangle}

\def\1{\hbox{I\hskip -2pt 1}}
\begin{document}
\maketitle
\begin{abstract}If $C_1$ is the convex hull of the curve of the  standard Brownian motion in the complex plane watched from $0$ to $1$, we consider the convex hulls 
of $C_1$ and several rotations of it and we compute the mean of the length of their perimeter by elementary calculations. 

\end{abstract}

\section{Introduction} If $C$ is a compact plane convex set, denote by $L(C)$ the length of its perimeter. 
Let $B(t)=(X(t),Y(t))$ be a standard Brownian motion in the plane and denote by $C_1(t)$ the convex hull of the curve $${\bf B}=\{B(s)\ ; 0\leq s\leq t\}.$$ This object has been considered by Paul L\'evy (1948), pages 239-240. The first author a long time ago has called attention on the fact that the mean of  $L(C_1(t))$ is $\sqrt{8\pi t}$ in Letac (1978). Here we consider the symmetric convex hull $C_2(t)$, namely the convex hull of $C_1(t)\cup -C_1(t)$. 
The convex set $C_2(t)$ should not be confused with the Minkowski sum $C_1(t)+(-C_1(t)),$ which is  bigger. 
By scaling we can guess the existence of a constant $\ell_2$  such that
$E[L(C_2(t))]=\ell_2\sqrt{8\pi t}.$ 
Similarly we consider  random convex sets obtained in the following way.
We let $\Omega\subset [0,2\pi[$ be a set of angles and
$C_\Omega$ the convex hull of $\cup_{\omega\in\Omega}R_\omega {\bf B}$, where $R_\omega$ is the rotation of angle $\omega$. Again by scaling there exists a constant $\ell_\Omega$  such that
$E[L(C_\Omega(t))]=\ell_\Omega\sqrt{8\pi t}.$ 
We shall compute these constants when
$\Omega$ is equal to one of the following sets
$$\{0\},\quad \{0,\pi\},\quad \{0,\pi/2\},\quad \{0,2\pi/3\},\quad \{0,\pi/2,\pi\}$$ or
$$\{0,2\pi/3,4\pi/3\}\quad
\{0,\pi/2,\pi,3\pi/2\},\quad [0,2\pi[$$
giving the result either in closed form or as an expression which can be evaluated numerically.
\section{Computations}
\subsection{Preliminary remarks}
Given a compact convex set $C$ in the Euclidean plane $\mathbb{R}^2$ consider its support function $h_C$ as the  function on $\mathbb{R}$ of period $2\pi$ defined by 
$$h_C(\theta)=\max_{(x,y)\in C}(x\cos\theta+y\sin\theta).$$
It is well known (see Letac (1983), formula (7.2) for a proof without regularity hypothesis) that the length of the perimeter of $C$ is $\int_0^{2\pi} h_C(\theta)d\theta.$ This formula is generally attributed to Cauchy. Let $\Omega\subset [0,2\pi[$ and define   $h_{{\bf B}}(\theta)=\max_{0\leq s\leq t}(X(s)\cos\theta+Y(s)\sin\theta)$ and 
$$h_\Omega(\theta)=\sup_{\omega\in\Omega}h_{{\bf B}}(\theta+\omega),$$ then rotation invariance of the law of Brownian motion implies
$$E[L(C_\Omega(t))]=2\pi E[h_\Omega(0)].$$
Let $P_{\Omega}$ be the convex set defined by the inequalities
$$x\cos\omega+y\sin\omega\leq 1;\omega\in \Omega.$$ In other terms, $P_{\Omega}$ is the intersection of half planes whose boundary is tangent to the unit circle.  Denote by 
 $T_\Omega$  the first time where Brownian motion exits this set.  Then a scaling argument shows that $h_\Omega(0)$ is distributed as $T_\Omega^{-1/2}$. It turns out that in several cases the Laplace transform 
$E[\exp(-\frac{\lambda^2}{2}T_\Omega)]$ is given by a simple formula. Under such a circumstance  
$$
E[h_\Omega(0)]=E[T_\Omega^{-1/2}]=\sqrt{\frac{2}{\pi}}\int_0^{\infty}E[\exp(-\frac{\lambda^2}{2}T_\Omega)]d\lambda
$$ 
and 
\begin{equation}\label{0}
E[L(C_\Omega)(t)]=\sqrt{8\pi t}\int_0^{\infty}E[\exp(-\frac{\lambda^2}{2}T_\Omega)]d\lambda.
\end{equation} \subsection{Some easy cases}
\subsubsection{$\Omega=\{0\}$}
The time $T_\Omega$ is the first hitting time of 1 for a linear Brownian motion, and
\begin{equation}
E[\exp(-\frac{\lambda^2}{2}T_\Omega)]=e^{-\lambda}.
\end{equation}
The Laplace tranform can be inverted to give
\begin{equation}\label{TAB}
P(T_\Omega\in dt)=\frac{e^{-1/2t}}{\sqrt{2\pi t^3}}.dt
\end{equation}
 Using (\ref{0}) we get
$E[L(C_1(t))]=\sqrt{8\pi t}.$ As expected: $l_{\Omega}=1.$
\subsubsection{$\Omega=\{0,\pi\}$} Now $T_\Omega$ is the exit time of linear Brownian motion from the interval $[-1,+1]$ and 
\begin{equation}\label{1}
E[\exp(-\frac{\lambda^2}{2}T_\Omega)]=1/\cosh\lambda
\end{equation} 
with density
\begin{equation}
P(T_\Omega\in dt)=\pi
\sum_{n=0}^\infty(-1)^n(n+\frac{1}{2})e^{-(n+\frac{1}{2})^2\pi^2t/2}dt
\end{equation}
and distribution function 
\begin{equation}\label{TAVB}
P(T_\Omega\geq  t)=\frac{2}{\pi}
\sum_{n=0}^\infty\frac{(-1)^n}{n+\frac{1}{2}}e^{-(n+\frac{1}{2})^2\pi^2t/2}
\end{equation}
see e.g. Biane, Pitman, Yor (2001).
Using again (\ref{0}) gives
$$E[L(C_2(t))]=\frac{\pi}{2}\sqrt{8\pi t}.$$
Thus $\ell_\Omega=\frac{\pi}{2}.$
\subsubsection{$\Omega=[0,2\pi[$}
The set $P_\Omega$ is the unit circle, and the
time $T_\Omega$ is the first hitting time of 1 by a Bessel process of dimension 2, therefore
\begin{equation}\label{F.-Kac}\mathbb{E}(e^{-\frac{\lambda^2}{2}T_1})=\frac{I_{0}(0)}{I_{0}(\lambda)}
\end{equation} where $$I_{0}(2z)=\sum_{k=0}^\infty\frac{z^{2k}}{k!^2}$$ is the Bessel function. Thus 
$$\ell_\Omega=\int_0^{\infty}\frac{d\lambda}{I_0(\lambda)}.$$
A  calculation by Mathematica (thanks to Daoud Bshouty for that) gives  
$\ell_\Omega=2,08323..$

\subsection{ A cone of angle $\pi/3$} We consider the case $\Omega=\{0,2\pi/3\}.$ Therefore
$T_\Omega$ is the exit time from a cone of angle $\pi/3$ tangent to the unit circle. In this case one can use results from Doumerc and O'Connell (2005). There, the exit time from the cone $x_1>x_2>x_3$ in ${\bf R}^3$, starting from the point $(x_1,x_2,x_3)$ is expressed as $$P(T\geq t)=p_{12}-p_{13}+p_{23}$$ with $p_{ij}=\sqrt{\frac{2}{\pi}}\int_0^{(x_i-x_j)/\sqrt{2t}}e^{-y^2/2}dy$.
Taking the orthogonal projection of Brownian motion on the hyperplane $x_1+x_2+x_3=0$,  with the starting point $(\sqrt2,0,-\sqrt 2)$, one gets
$$P(T_\Omega\geq t)=\sqrt{\frac{2}{\pi}}\left (2\int_0^{1/\sqrt t}e^{-y^2/2}dy-
\int_0^{2/\sqrt t}e^{-y^2/2}dy\right)$$ and a straightforward computation gives
$$E[T_\Omega^{-1/2}]=\frac{3}{2}\sqrt{\frac{2}{\pi}}$$
and
$$E[L(C_\Omega(t))]=\frac{3}{2}\sqrt{8\pi t}.$$ Thus $\ell_\Omega=3/2.$

\subsection{An equilateral triangle}
An interesting case is when $\Omega=\{0,2\pi/3,4\pi/3\}$ and 
$P_\Omega$ is an equilateral triangle centered at 0. We will compute the distribution of $T_\Omega$ and get a curious identity in law.

In order to do the computation it will be convenient to  take the equilateral triangle $\Delta$  with vertices $0,1,-j^2$ where $j=(-1+i\sqrt{3})/2$. Its   center is  $x_0=(1-j^2)/3$. The case of the triangle $P_\Omega$ can be obtained by an affinity.
The orthogonal reflections with respect to the three lines bordering the triangle
$\Delta$ generate a group $W$ of affine isometries of the Euclidian plane, whose fundamental domain is $\Delta$. Using the reflection principle,
we can compute the probability transition of Brownian motion in the triangle, killed at the boundary. One obtains
$$p^0_t(x,y)=\sum_{w\in W}\det(w)p_t(x,w(y))$$
where $p_t(x,y) $ is the ordinary heat kernel on the Euclidean plane and $x$ and $y$ are points of this plane. .
If we denote $R$ the lattice generated by $1$ and $j$, then the dual lattice $\hat R$ (such that
$\langle w,\hat w\rangle\in{\bf Z}$) is generated by the dual basis
$a=1+\frac{i}{\sqrt{3}}$ and $b=\frac{2i}{\sqrt{3}}$.
The plane is tiled by translates of $\Delta$ by $R$, which correspond to direct isometries   $w\in W$ and translates of $\bar\Delta$.
Thus 
$$p^0_t(x,y)=\sum_{r\in R}p_t(x-y,r)-p_t(x-\bar y,r).$$
The function 
$$h(t,z)=\sum_{r\in R}p_t(z,r)$$ is periodic in $z$ with period lattice $R$ and $p^0_t(x,y)=h(t,x-y)-h(t,x-\bar y).$ 
One can apply Poisson summation formula, i.e. expand it in terms of the 
characters $e^{2i\pi\langle .,\hat r\rangle}$
obtaining
$$h(t,z)=\frac{2}{\sqrt{3}}
\sum_{\hat r\in \hat R}e^{2i\pi\langle z,\hat r\rangle}
e^{-2\pi^2|\hat r|^2t}$$
Let us now integrate with respect to $y$ in the triangle to obtain the distribution function for $T_\Omega$, the first exit time of Brownian motion starting from $x_0$,
$$P(T_\Omega>t)=\int_{\Delta}p_t^0(x_0,y)dy.$$

We use the coordinates
$y=y_1+y_2j,1\geq  y_1\geq y_2\geq 0$ on the triangle, with Jacobian
$\frac{\sqrt{3}}{2} $  and 
$\hat r=ma+nb$ on $\hat R$, then $$\langle y,\hat r\rangle=my_1+ny_2,\qquad 
\langle x_0,\hat r\rangle=2m/3+n/3,\qquad
|\hat w|^2=4(m^2+mn+n^2)/3.$$ Furthermore, by periodicity
$$\int_\Delta h(t,x-\bar y)dy=\int_\Delta h(t,x-j-\bar y)dy=\int_{\bar\Delta+j} h(t,x- y)dy$$ and the triangle
$\bar \Delta+j$ corresponds to the coordinates
$1\geq y_2\geq y_1\geq 0$. Let us compute
\begin{eqnarray*}
\int_\Delta h(t,x_0- y)dy&=&\frac{\sqrt{3}}{2}\sum_{n,m}e^{2i\pi (2m/3+n/3)}
\int_0^1dy_1\int_0^{y_1}e^{-2i\pi(my_1+ny_2)}
dy_2e^{-8\pi^2(m^2+mn+n^2)t/3}
\end{eqnarray*}
It is easy to see that terms with $n\ne 0,m\ne 0,m+n\ne 0$ give $0$ contribution. Also the term
$m=n=0$ will be cancelled by the other integral. The  sequence
$m=0,n\ne 0$,
gives 
\begin{eqnarray*}&&\sum_{n\ne 0}e^{2i\pi n/3}
\int_0^1dy_1\int_0^{y_1}e^{-2i\pi ny_2}
dy_2e^{-8\pi^2n^2t/3}\\
&=&\sum_{n\ne 0}e^{2i\pi n/3}
\int_0^1\frac{1}{2i\pi n}(1-e^{-2i\pi ny_1})dy_1
e^{-8\pi^2n^2t/3}\\&=&
\sum_{n\ne 0}e^{2i\pi n/3}
\frac{1}{2i\pi n}
e^{-8\pi^2n^2t/3}
\end{eqnarray*}

The  sequence
$m\ne 0,n=0$,
gives 
\begin{eqnarray*}&&\sum_{m\ne 0}e^{4i\pi m/3}
\int_0^1dy_1e^{-2i\pi my_1}\int_0^{y_1}
dy_2e^{-8\pi^2m^2t/3}\\
&=&\sum_{m\ne 0}e^{4i\pi m/3}
\int_0^1dy_1y_1e^{-2i\pi y_1m}
e^{-8\pi^2m^2t/3}\\&=&
-\sum_{m\ne 0}e^{4i\pi m/3}
\frac{1}{2i\pi m}
e^{-8\pi^2m^2t/3}
\end{eqnarray*}
and the  sequence
$m\ne 0,m+n=0$,
gives 
\begin{eqnarray*}&&\sum_{m\ne 0}e^{2i\pi m/3}
\int_0^1dy_1\int_0^{y_1}e^{2i\pi(-my_1+my_2)}
dy_2e^{-8\pi^2m^2t/3}\\
&=&\sum_{m\ne 0}e^{2i\pi m/3}
\frac{1}{2i\pi m}
\int_0^1dy_1(1-e^{-2i\pi y_1m})
e^{-8\pi^2m^2t/3}\\&=&
\sum_{m\ne 0}e^{2i\pi m/3}
\frac{1}{2i\pi m}
e^{-8\pi^2m^2t/3}
\end{eqnarray*}
The sum of these three terms is
$$3\sum_{n=1}^\infty \frac{\sin(2\pi n/3)}{\pi n}e^{-8\pi^2n^2t/3}.$$
 The other integral is
$$\int_\Delta h(t,x_0- \bar y)dy=\sum_{n,m}e^{2i\pi (2m/3+n/3)}
\int_0^1dy_2\int_0^{y_2}e^{-2i\pi(my_1+ny_2)}
dy_1e^{-8\pi^2(m^2+mn+n^2)t/3}$$
Using 
$$
\int_0^1dy_1\int_0^{1}e^{-2i\pi(my_1+ny_2)}
dy_1dy_2=0$$
if $m$ or $n$ is $\ne0$ we see that, except for the term $m=n=0$ the other terms are the opposite of what we computed, therefore after some rearrangement

$$P(T_\Omega>t)=3\sqrt{3}\sum_{n=1}^{\infty}\frac{\chi_3(n)}{\pi n}e^{-8\pi^2n^2t/3}$$
where $$\chi_3(n)=\frac{2}{\sqrt{3}}\sin(2\pi n/3)$$ is the multiplicative Dirichlet character modulo 3.
 The density is
$$P(T_\Omega\in dt)=8\sqrt{3}\sum_{n=1}^{\infty}\pi n\chi_3(n)e^{-8\pi^2n^2t/3}.$$
\subsubsection{A curious identity in law}
Let  $S$ be the first hitting time of $3a$ by a three-dimensional
 Bessel process starting from $a$. The  Laplace transform is
$$E[e^{-\lambda^2S/2}]=\frac{3\sinh(\lambda a)}{\sinh(3\lambda a)}
=3\sum_{n=1}^{\infty}\chi_3(n)e^{-2n\lambda a}.$$
Inverting term by term gives the density
$$\sum_{n=0}^{\infty}6an\frac{\chi_3(n)e^{-2n^2a^2/t}}
{\sqrt{2\pi t^3}}$$
We can also use the alternative expression
$$\frac{3\sinh(\lambda a)}{\sinh(3\lambda a)}=\sum_{n=1}^\infty\frac{3\sqrt{3}\pi n\chi_3(n)}{\pi^2n^2+9\lambda^2a^2}$$
and invert term by term the Laplace transform to get
$$P(S\in dt)=\sum_{n=1}^{\infty}\frac{3\sqrt{3}\pi n\chi_3(n)}{18a^2}
e^{-\pi^2n^2t/(18a^2)}$$
The agreement between these two expressions is an instance of the functional equation of theta series.
Putting $8/3=1/(18a^2)$, or $a=\frac{1}{4\sqrt{3}}$ we see that the exit time $T_\Omega$ and the time $S$ have the same distribution. Furthermore, the Mellin transform of this distribution is
$$E(T_\Omega^s)=\pi^{-2s}8^{1-s}3^{s+1/2}\Gamma(s+1)L_{\chi_3}(2s+1)$$
where $L_{\chi_3}$ is the Dirichlet L-function associated with the character $\chi_3$.
Coming back to $T_\Omega$, the triangle $P_\Omega$ has sides of size $2\sqrt{3}$ therefore 
$$E[e^{-\frac{\lambda^2}{2}T_\Omega}]=\frac{3\sinh(\lambda/2)}{\sinh(3\lambda/2 )}.
$$
Using (\ref{0}) we are reduced to an elementary integral, and the result is
$$\ell_\Omega=\frac{\pi}{\sqrt{3}}.$$

\subsection{The other cases}
\subsubsection{$\Omega=\{0,\pi/2\}$} 
We have 
$$h_{\Omega}(\theta)=\max\{X(s)\cos \theta+Y(s)\sin\theta,\ Y(s)\cos \theta-X(s)\sin\theta\ ; \ 0\leq s\leq t\}.$$ The trick is to use the fact that $(X(s)\cos \theta+Y(s)\sin\theta)_{s\geq 0}$  and $(Y(s)\cos \theta-X(s)\sin\theta)_{s\geq 0}$  are two   standard one dimensional Brownian motions
which are independent. As a consequence, using (\ref{TAB}) we have 
$$\Pr(h_{\Omega}(\theta)\leq h)=H^2(h/\sqrt{t})$$
where  
\begin{equation}\label{A1}H(z)=-1+2\int_{-\infty}^ze^{-\frac{x^2}{2}}\frac{dx}{\sqrt{2\pi}}.\end{equation}
 It leads to $$\mathbb{E}(L(\Omega))=2\pi\int_0^{\infty}(1-H^2(h/\sqrt{t}))dh=2\pi\sqrt{t}\int_0^{\infty}(1-H^2(z))dz.$$We write 
$$1-H^2(z)=4\int_z^{\infty}e^{-\frac{x^2}{2}}\frac{dx}{\sqrt{2\pi}}\left(1-\int_z^{\infty}e^{-\frac{y^2}{2}}\frac{dy}{\sqrt{2\pi}}\right)$$
By the change of variables $x=zu$ and $y=zv$ and $w=z^2/2$ we get easily  
$$\mathbb{E}(L(K_2(t)))=\sqrt{8\pi t}\left(2-\int_1^{\infty}\int_1^{\infty}\frac{dudv}{(u^2+v^2)^{3/2}}\right)$$
Introduce the function $(x,y)\mapsto g(x,y)$ defined for $x$ and $y$ in $\mathbb{R}\setminus \{0\}$ by \begin{equation}\label{G}g(x,y)=\frac{\sqrt{x^2+y^2}}{xy}.\end{equation}  A tedious calculation shows that if $a<b$ and $c<d$ with $(0,0)\notin [a,b]\times [c,d]$  and if furthermore $a,b,c,d\neq 0$ we have
\begin{equation}\label{H}\int_{a}^{b}\int_{c}^{d}\frac{dudv}{(u^2+v^2)^{3/2}}=g(a,c)+g(b,d)-g(a,d)-g(b,c).\end{equation}
From this we get that $\int_1^{\infty}\int_1^{\infty}\frac{dudv}{(u^2+v^2)^{3/2}}=2-\sqrt{2}$ which finally gives that $$\ell_\Omega=\sqrt{2}.$$

\subsubsection{$\Omega=\{0,\pi/2,\pi,3\pi/2\}$} The method is quite similar to the method used before, but the distribution function $H$ of the maximum of a Brownian motion is replaced by the distribution function $L$ of the maximum of its absolute value, which can be obtained from (\ref{TAVB}). Still with  the previous notations we have 
$$h_{C_\Omega}(\theta)=\max\{|X(s)\cos \theta+Y(s)\sin\theta|,\ |Y(s)\cos \theta-X(s)\sin\theta|\ ; \ 0\leq s\leq t\}.$$ We use again the fact that $(X(s)\cos \theta+Y(s)\sin\theta)_{s\geq 0}$  and $(Y(s)\cos \theta-X(s)\sin\theta)_{s\geq 0}$  are two   standard one dimensional independent Brownian motions. We get
$$\Pr(h_{C_\Omega}(\theta)\leq h)=L^2(h/\sqrt{t})$$ and  $$\mathbb{E}(L(C_\Omega(t)))=\int_0^{2\pi}\mathbb{E}(h_{C_4(t)}(\theta))d\theta=2\pi\int_0^{\infty}(1-L^2(h/\sqrt{t}))dh=2\pi\sqrt{t}\int_0^{\infty}(1-L^2(z))dz.$$ We now  write
$$1-L^2(z)=(1-L(z))(1+L(z))=\frac{2}{\pi}\sum_{k=-\infty}^{\infty}\sum_{n=-\infty}^{\infty}\int_{(4k-1)z}^{(4k+1)z}e^{-\frac{x^2}{2}}dx\int_{(4n+1)z}^{(4n+3)z}e^{-\frac{y^2}{2}}dy.$$ We are led to the computation of the integral for fixed $(k,n)\in \mathbb{Z}^2$
$$I(k,n)=\sqrt{\frac{2}{\pi}}\int_0^{\infty}\left(\int_{(4k-1)z}^{(4k+1)z}e^{-\frac{x^2}{2}}dx\int_{(4n+1)z}^{(4n+3)z}e^{-\frac{y^2}{2}}dy\right)dz.$$
By the same change of variables as before $x=zu$ and $y=zv$ and $w=z^2/2$ we get   
\begin{equation}\label{C}I(k,n)=\int_{4k-1}^{4k+1}\int_{4n+1}^{4n+3}\frac{dudv}{(u^2+v^2)^{3/2}}.\end{equation}
In order to compute $I(k,n)$ we use the function $g$ defined  by (\ref{G}) and the identity (\ref{H}).  
From (\ref{C}) it is apparent that $I(k,n)>0$. To see also that \begin{equation}\label{E}\ell_\Omega=\sum_{k=-\infty}^{\infty}\sum_{n=-\infty}^{\infty}I(k,n)<\infty\end{equation}
we use (\ref{C}) again for observing that $\ell_\Omega$ is the integral of the function $(u^2+v^2)^{-3/2}$ over the union of squares 
$$D=\cup_{k=-\infty}^{\infty}\cup_{n=-\infty}^{\infty}[4k-1,4k+1]\times [4n+1,4n+3]$$ which satisfies $D\subset A=\{(u,v)\ ;\ u^2+v^2\geq 1\}.$ This implies, using polar coordinates $u=r\cos \theta$ and $v=r\sin \theta:$
$$\ell_\Omega=\int_D\frac{dudv}{(u^2+v^2)^{3/2}}< \int_A\frac{dudv}{(u^2+v^2)^{3/2}}=2\pi\int_1^{\infty}r^{-2}dr=2\pi.$$  In order to have some numerical evaluation of $\ell_\Omega$ one can observe that we have the following two symmetries: for all $(k,n)\in\mathbb{Z}^2$
$$I(k,n)=I(-k,n)=I(k,-n-1).$$
This implies that $\ell_\Omega$ can be written as $\ell_\Omega=S_0+S_1$ with
\begin{equation}\label{J}S_0=\sum_{n=-\infty}^{\infty}I(0,n)=4\sum_{n=0}^{\infty}\left( \frac{\sqrt{1+(4n+1)^2}}{4n+1} - \frac{\sqrt{1+(4n+3)^2}}{4n+3}\right)=4\sum_{m=0}^{\infty}(-1)^m u_m\end{equation}
with $u_m=\sqrt{1+\frac{1}{(2m+1)^2}} - 1,$ and 
$S_1=4\sum_{n=0}^{\infty}\sum_{k=1}^{\infty}I(k,n).$
Since $(u_m)_{m\geq 0}$ is decreasing we can estimate $\frac{S_0}{4}\approx0.3725...$
We compute $\ell_\Omega$ by Mathematica (we thank Daoud Bshouty for this) and we get $$\ell_\Omega=1,9178...$$ 
\subsubsection {$\Omega=\{0,\pi/2,\pi\}$}
Again the same method gives
$$\mathbb{E}(L(C_\Omega))=2\pi\int_0^{\infty}(1-L(h/\sqrt{t})H(h/\sqrt{t}))dh=2\pi\sqrt{t}\int_0^{\infty}(1-L(z)H(z))dz.$$
We write for simplicity $L=1-2L_1$ and $H=1-2H_1$ so that 
$$\mathbb{E}(L(C_\Omega))=2\pi\sqrt{t}\int_0^{\infty}\left(2L_1(z)+2H_1(z)-4L_1(z)H_1(z)\right)dz.$$ A previous calculation  has shown that 
$2\pi\sqrt{t}\int_0^{\infty}2L_1(z)dz=\frac{\pi}{2}\sqrt{8\pi t}$ and Letac (1978) -or a direct calculation- shows $2\pi\sqrt{t}\int_0^{\infty} 2H_1(z)dz=\sqrt{8\pi t}.$ The delicate remaining point is the computation of $\int_0^{\infty}4L_1(z)H_1(z)dz.$ By the same methods as before,  the same changes of variable $x=zu$ and $y=zv$ and $w=z^2/2$, the use of the function $g$ of (\ref{G}) and the formula  (\ref{H}) we get 

\begin{eqnarray*}2\pi\sqrt{t}\int_0^{\infty}4L_1(z)H_1(z)dz&=&
\sqrt{8\pi t}\sum_{n=-\infty}^{\infty}\int_1^{\infty}\int_{4n+1}^{4n+3}\frac{dudv}{(u^2+v^2)^{3/2}}\\&=&\sqrt{8\pi t}\sum_{n=-\infty}^{\infty}(g(1,4n+1)-g(1,4n+3)).
\end{eqnarray*}
It is easily seen that with the notation (\ref{J})
\begin{equation}\label{F}\sum_{n=-\infty}^{\infty}(g(1,4n+1)-g(1,4n+3))=\frac{S_0}{2}\approx 0.7450...\end{equation}
Thus $$\ell_\Omega\approx\frac{\pi}{2}+0.2550....$$

\section{Comments}

\noindent1. The explicit calculations of the paper for the mean perimeter are feasible since we are able to have the exact distribution of the  support functions $\theta\mapsto h_C(\theta)$ for the random convex sets considered here. The general case seems a difficult problem. 

\noindent2. Calculations for the mean areas rely on the formula $A(C)=\int_{0}^{2\pi}(h_C(\theta)^2-h_{C}'(\theta)^2)d\theta$, where $h_{C}'(\theta)$ denotes the derivative on the left, which can be shown to exist (see Letac (1983)) and   ultimately of the distribution of $h'_C(\theta)$ when $C$ is random.  This problem was solved in the unpublished thesis of El Bachir (1983) for $C_1(t)$, where $\mathbb{E}(A(C_1(t))=\pi t/2$ is proved.  It is unsolved for the other convex sets considered here, except for $C_{\infty}(t)$ the smallest circle, centered at 0,   surrounding $C_1(t)$, indeed
$$\mathbb{E}(A(C{_\infty}(t))=\pi t \int_0^{\infty}\frac{ds}{I_0(\sqrt{2s})}=\pi t\times 3.06883.$$ 
a little more than six times the quantity $\mathbb{E}(A(C{_1}(t)).$ Thanks to Daoud Bshouty for the numerical evaluation of this integral.

\section{References}
\vspace{4mm}\noindent\textsc{Biane, P.; Pitman, J.; Yor, M. } (2001) 'Probability laws related to the Jacobi theta and Riemann zeta functions, and Brownian excursions' \textit{Bull. Amer. Math. Soc. (N.S.)}
 \textbf{38}, no. 4, 435--465.

    \vspace{4mm}\noindent\textsc{Doumerc, Y.: O'Connell, N } (2001) 'Exit problems associated with finite reflection groups' \textit{Probab. Th. Rel. Fields)}
 \textbf{132}, no. 4, 501--538.

\vspace{4mm}\noindent
\textsc{El Bachir, M.} (1983) \textit{L'enveloppe convexe du mouvement brownien.} Th\`ese de 3\`eme cycle, Universit\'e Paul Sabatier, Toulouse.

 \vspace{4mm}\noindent\textsc{Feller, W. } (1966) \textit{An Introduction to Probability Theory and Its Applications}, Vol 2. Wiley, New York.

\vspace{4mm}\noindent \textsc{Letac, G.} (1978)
 'Advanced problem 6230'
 {\it Am. Math. Monthly}.  Solution by L. Tak\'acs {\bf
87},  142 (1980).

 \vspace{4mm}\noindent\textsc{Letac, G.} (1983) 'Mesures sur le cercle et convexes du plan' \textit{Annales scientifiques de l'Universit\'e de Clermont-Ferrand II}
 \textbf{76}, 35-65.

 \vspace{4mm}\noindent\textsc{L\'evy, P.} (1948) \textit{Processus stochastiques et mouvement brownien}, Gauthier Villars, Paris.

\end{document}